\documentclass[12pt]{amsart}

\usepackage{tikz,url}
\usepackage{amsmath,amssymb}

\usepackage[margin = 2.5cm]{geometry}

\theoremstyle{theorem}
\newtheorem{theorem}{Theorem}
\newtheorem{observation}[theorem]{Observation}

\theoremstyle{definition}

\allowdisplaybreaks

\makeatletter
\@addtoreset{footnote}{page}
\makeatother

\DeclareMathOperator{\des}{des}
\DeclareMathOperator{\NN}{\mathbb{N}}
\DeclareMathOperator{\A}{\mathcal{A}}
\DeclareMathOperator{\LL}{\mathcal{L}}
\DeclareMathOperator{\op}{\Omega}

\newcommand{\el}[2]{\genfrac{\langle}{\rangle}{0pt}{}{#1}{#2}}
\newcommand{\st}[2]{\genfrac{\{}{\}}{0pt}{}{#1}{#2}}
\newcommand{\fl}{\begin{tikzpicture}[baseline=0,scale=0.25] \draw[line width=1,cap=round] (0,0) --(0,1) node[circle,fill=black,inner sep=1]{}; \end{tikzpicture}}
\newcommand{\anc}{\begin{tikzpicture}[baseline=1pt,scale=0.25] \draw[line width=1,cap=round] (0,1) --(0,0) node[circle,fill=black,inner sep=1]{}; \end{tikzpicture}}

\begin{document}

\title{The Eulerian numbers can D.I.E.}

\author{Matja\v z Konvalinka and T. Kyle Petersen}

\begin{abstract}
The Eulerian numbers form a triangular array with many interesting properties. The numbers arise from various combinatorial and probabilistic interpretations, and have been studied in a variety of mathematical contexts. In this article we examine two distinct alternating sign formulas for the Eulerian numbers and show how they can be proved using a sign-reversing involution technique described by Benjamin and Quinn known as the ``D.I.E.'' method. Each of these arguments lends itself to a broad generalization, shedding light on different parts of mathematics.
\end{abstract}

\maketitle

The \emph{Eulerian numbers} form a triangular array similar to Pascal's triangle with many pleasing features. The numbers arise from various combinatorial, algebraic, and probabilistic interpretations, and have been studied in a variety of mathematical contexts. See \cite{Petersen} for an introduction to the subject. The numbers are defined combinatorially as 
\[
 \el{n}{k} = |\{ \pi \in S_n : \des(\pi) = k\}|,
\]
where $S_n$ denotes the set of permutations of $\{1,2,\ldots,n\}$ and $\des(\pi)$ is the number of \emph{descents} of $\pi$, i.e., $\des(\pi) = |\{ i : \pi(i) > \pi(i+1) \}|$. For example, if $n=3$, there are six permutations, written here in one-line notation $\pi(1)\pi(2)\pi(3)$:
\[
 123, 132, 213, 231, 312, 321,
\]
and they have $0, 1, 1, 1, 1,$ and $2$ descents, respectively. Thus, $\el{3}{0} = 1$, $\el{3}{1} = 4$, and $\el{3}{2}=1$. See Table \ref{tab:Eul} for the values of $\el{n}{k}$ with $1\leq k+1 \leq n \leq 8$.

Given that one of us has written a book about the Eulerian numbers \cite{Petersen}, the title of this article may seem a little surprising. But the title does not indicate any ill feeling toward this famous array of numbers---to the contrary! 
Inspired by a question from the first author we have found a beautiful proof of the following identity
\begin{equation}\label{eq1}
\el{n}{k} = \sum_{i=0}^k (-1)^i\binom{n+1}{i}(k+1-i)^n.
\end{equation}
The first author had initially hoped to use the method of ``sijections'' (signed bijections), as developed in recent work on the Alternating Sign Matrix theorem \cite{FK1, FK2}, to prove Equation \eqref{eq1}. In fact the sijection approach \emph{can} be applied to this formula, but in our opinion it turned out to be less appealing than the proof provided in the next section of this paper. There we use the D.I.E. method described by Benjamin and Quinn \cite{BQ}, which we state almost verbatim here:
\begin{itemize}
\item \textbf{Describe:} Describe a set of objects that is being counted by the sum when we ignore the sign.
\item \textbf{Involution:} Find an involution on the set of objects that matches items that are counted positively and items that are counted negatively.
\item \textbf{Exception:} Describe the exceptions that are mapped to themselves. Count these exceptions and note their sign.
\end{itemize}

After this initial success, we go further in the subsequent section to describe a D.I.E. proof for another alternating sign identity for Eulerian numbers. This one relates Eulerian numbers, binomial coefficients, and the \emph{Stirling numbers of the second kind}, denoted $\st{n}{k}$. Here $\st{n}{k}$ counts the number of set partitions of an $n$-element set into $k$ (unordered) parts. See Table \ref{tab:Stir} for the values of $\st{n}{k}$ with $1\leq k \leq n \leq 8$. The identity we will prove is
\begin{equation}\label{eq2}
 \el{n}{k} = \sum_{i=1}^{k+1} (-1)^{k+1-i}\binom{n-i}{k+1-i}\st{n}{i}i!.
\end{equation}

\begin{table}[h]
\caption{The Eulerian numbers $\el{n}{k}$, $1\leq k+1\leq n \leq 8$.}\label{tab:Eul}
\begin{center}
\begin{tabular}{r | c c c c c c c c c c c c}
\hline
$n \backslash k$ & 0 & 1 & 2 & 3 & 4 & 5 & 6 & 7\\
\hline
1 & 1 &&&&& \\
2 & 1 & 1 &&&&&\\
3 & 1 & 4 & 1 &&\\
4 & 1 &  11 & 11 & 1 \\
5 & 1 & 26 & 66 & 26 & 1\\
6 & 1 & 57 & 302 & 302 & 57 & 1\\
7 & 1 & 120 & 1191 & 2416 & 1191 & 120 & 1\\
8 & 1 & 247 & 4293 & 15619 & 15619 & 4293 & 247 & 1\\
\hline
\end{tabular}
\end{center}
\end{table}

\begin{table}[h]
\caption{The Stirling numbers of the second kind $\st{n}{k}$, $1\leq k \leq n \leq 8$.}\label{tab:Stir}
\begin{center}
\begin{tabular}{r | c c c c c c c c c c c c}
\hline 
$n \backslash k$ & 1 & 2 & 3& 4 & 5 & 6 & 7 & 8\\
\hline
1 & 1 &&&&& \\
2 & 1 & 1 &&&&&\\
3 & 1 & 3 & 1 &&\\
4 & 1 & 7 & 6 & 1 \\
5 & 1 & 15 & 25 & 10 & 1\\
6 & 1 & 31 & 90 & 65 & 15 & 1\\
7 & 1 & 63 & 301 & 350 & 140 & 21 & 1\\
8 & 1 & 127 & 966 & 1701 & 1050 & 266 & 28 & 1\\
\hline
\end{tabular}
\end{center}
\end{table}

We remark that both identities are well established in the literature, e.g., \eqref{eq1} and \eqref{eq2} are found as Equations (6.38) and (6.40), respectively, on page 269 of \cite{GKP}.\footnote{We note that \eqref{eq2} differs from (6.40) only by the symmetry $\el{n}{k}=\el{n}{n-1-k}$ and re-indexing.} They are usually proved by a combination of techniques, e.g., \eqref{eq1} follows by inverting an identity known as ``Worpitzky's'' or ``Li Shan-Lan's'' identity:
\[
 (k+1)^n = \sum_{i=0}^{n-1} \el{n}{i}\binom{k+n-i}{n},
\]
which itself can be proved bijectively, inductively, or with generating function techniques. Equation \eqref{eq2} follows from a similar identity connecting the Eulerian numbers to the Stirling numbers of the second kind:
\[
 k!\st{n}{k} = \sum_{i=0}^{k-1} \el{n}{i} \binom{n-1-i}{k-1-i},
\]
which can be proved bijectively or by comparing the generating function for the Eulerian numbers to the generating function for the Stirling numbers of the second kind. 

Our contribution in this article is to provide what we understand to be new proofs of these identities using the appealing combinatorial paradigm suggested by Benjamin and Quinn. Further, both identities are special cases of more general identities that admit proofs using the same D.I.E. approach. 

In the case of Equation \eqref{eq1}, the more general identity can be stated in the language of Stanley's theory of $P$-partitions \cite[Section 3.15]{Stanley}, and has a geometric interpretation as counting lattice points in dilations of \emph{order polytopes}. Here we get an alternating sign formula for \emph{$P$-Eulerian numbers} $\el{P}{k}$, where $P$ is a finite partially ordered set. 

For Equation \eqref{eq2}, the more general identity is stated in the language of combinatorial topology, as a statement about counting faces in a \emph{partitionable simplicial complex}. This result requires some familiarity with abstract simplicial complexes, and we introduce the reader to the relevant definitions. After establishing the general case, we show how considering the barycentric subdivision of a simplex yields a third alternating sign formula for Eulerian numbers:
\begin{equation}\label{eq3}
 \el{n}{k} = \sum_{i=0}^k (-1)^{k-i} \binom{n-i}{k-i}\st{n+1}{i+1}i!.
\end{equation}

For an example of each alternating sign identity, we compute $\el{5}{2} = 66$ in three different ways. From \eqref{eq1} we find:
\[
 \el{5}{2} = \binom{6}{0}3^5 - \binom{6}{1}2^5 + \binom{6}{2}1^5 = 243-192+15 = 66,
\]
from \eqref{eq2} we find:
\[
 \el{5}{2} = \binom{4}{2}\st{5}{1}1! - \binom{3}{1}\st{5}{2}2! + \binom{2}{0}\st{5}{3}3! = 6 - 90 + 150 = 66,
\]
and from \eqref{eq3} we find:
\[
 \el{5}{2} = \binom{5}{2}\st{6}{1}0! - \binom{4}{1}\st{6}{2}1!+\binom{3}{0}\st{6}{3}2! = 10 - 124 + 180 = 66.
\]

We hope this article gives more evidence for the usefulness of the D.I.E. method. We also hope that it showcases two important contexts in which the Eulerian numbers arise. As suggested in \cite{Petersen}, the Eulerian numbers can make for an accessible entry point to many interesting areas of mathematics!

\section{$P$-partitions and geometry}

This first main section discusses Equation \eqref{eq1} and its generalization to the realm of partially ordered sets. We now proceed to our first proof to D.I.E. for.

\subsection{First D.I.E. proof}

We begin by describing the set of \emph{barred permutations}, which is equivalent to the set of placements of balls into boxes. Suppose we have $k+1$ boxes, labeled $0,1,\ldots,k$, and we put $n$ labeled balls into these boxes. For example, with $n=9$ and $k=7$, we might have
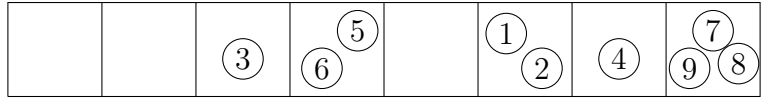
\[
 \begin{tikzpicture}[scale=1.25]
  \draw (0,0) -- (0,1) -- (8,1)--(8,0)--(0,0);
  \foreach \x in {1,...,7}{
   \draw (\x,0)--(\x,1);
  }
  \draw (2.5,0.4) node[circle,draw=black,inner sep=2] {3};
  \draw (3.34,0.3) node[circle,draw=black,inner sep=2] {6};
  \draw (3.72,0.7) node[circle,draw=black,inner sep=2] {5};
  \draw (5.3,0.67) node[circle,draw=black,inner sep=2] {1};
  \draw (5.68,0.3) node[circle,draw=black,inner sep=2] {2};
  \draw (6.5,0.4) node[circle,draw=black,inner sep=2] {4};
  \draw (7.5,0.7) node[circle,draw=black,inner sep=2] {7};
  \draw (7.77,0.35) node[circle,draw=black,inner sep=2] {8};
  \draw (7.25,0.3) node[circle,draw=black,inner sep=2] {9};
 \end{tikzpicture}
\]
which we can write more succinctly as the barred permutation 
\[
 \beta=||3|56||12|4|789.
\]
Notice that we list elements from each box in increasing order, and when there are $k+1$ boxes, there are $k$ bars in $\beta$. Let $\pi = \pi(\beta)$ in $S_n$ denote the \emph{underlying permutation} for the barred permutation $\beta$ of size $n$. Let $B_{n,k}$ denote the set of barred permutations of $n$ elements with $k$ bars.


There are $(k+1)^n$ ways to put $n$ labeled balls into $k+1$ labeled boxes, so we make the following observation.

\begin{observation}\label{obs1}
Fix $n\geq 1$ and $k\geq 0$. Then $|B_{n,k}| = (k+1)^n$.
\end{observation}

Because elements between bars are listed in increasing order, we can see that descents of $\pi$ can only occur between elements from different boxes, i.e., there must be at least one bar in each descent position of $\pi$. We make the following observation.

\begin{observation}\label{obs2}
Suppose $\beta \in B_{n,k}$. Then $\des(\pi(\beta))\leq k$.
\end{observation}

Now define the set of \emph{anchored} barred permutations. To explain this, we make three types of bars. There are:
\begin{itemize}
\item \emph{anchor bars}, indicated with $\anc$, which appear exactly once in gaps with a descent and nowhere else,
\item \emph{unnecessary bars}, indicated with $|$, which can appear in any gap any number of times, and
\item \emph{float bars}, indicated with $\fl$, which can appear at most once per gap.
\end{itemize}
If bars of different types appear in the same gap, we place float bars to the left of unnecessary bars, which are to the left of anchor bars.

For example,
\[
\beta=\fl |3|56|\anc 12\fl 4|789
\]
is an anchored barred permutation with underlying permutation $\pi=356124789$. It has $7$ bars: one anchor, four unnecessary bars, and two floats. Notice that the number of anchors is precisely the number of descents in the underlying permutation, and thus we make the following observation.

\begin{observation}\label{obs3}
The number of anchored barred permutations of size $n$ with $k$ anchors and no other bars is $\el{n}{k}$.
\end{observation}

Let $B_{n,k}^i$ denote the set of anchored barred permutations, $i$ of which are floats. The \emph{sign} of an anchored barred permutation $\beta$ is $(-1)^{\# \textrm{float bars}}$. Ignoring the floats, we see there are $k-i$ other bars in such a barred permutation. From Observation \ref{obs1} we know there are $(k+1-i)^n$ barred permutations with $k-i$ bars, and we can choose to add one float bar to any $i$ of the $n+1$ gaps. Our next observation summarizes our enumeration of the elements of $B_{n,k}^i$. 

\begin{observation}\label{obs4}
Fix $n\geq 1$, $k\geq 0$ and $i\leq k$. Then $|B_{n,k}^i| = \binom{n+1}{i}(k+1-i)^n$. Moreover, each $\beta$ in $B_{n,k}^i$ has sign $(-1)^i$.
\end{observation}

Thus our \textbf{description} for proving Equation \eqref{eq1} is to define the set $\mathcal{B}_{n,k}$ of all anchored barred permutations of size $n$ with $k$ bars, i.e.,
\[
 \mathcal{B}_{n,k} = \bigcup_{i=0}^{k} B_{n,k}^i.
\]

Our involution will toggle between whether the leftmost non-anchor bar is a float or an unnecessary bar: $\fl{} \leftrightarrow |$, resulting in a new anchored barred permutation $\beta'$. Here are three examples of $\beta \leftrightarrow \beta'$, all with $n=9$ and $k=7$:
\[
\begin{array}{rcl}
\\
\fl |3|56|\anc 12\fl 4|789 & \leftrightarrow & ||3|56|\anc 12\fl 4|789 \\
\\
 3 || 56\anc 1 \fl 24\fl 789 || & \leftrightarrow & 3 \fl | 56\anc 1 \fl 24\fl 789 || \\
 \\
 79\anc 3 \anc 24|\anc 19 |\anc 6\anc 5 & \leftrightarrow & 79\anc 3 \anc 24 \fl \anc 19 |\anc 6\anc 5 \\
 {}\\
\end{array}
\]
Notice that by construction, $\beta$ and $\beta'$ have opposite sign.

The \textbf{involution} $\iota_1: \mathcal{B}_{n,k} \to \mathcal{B}_{n,k}$ is as follows:
\[
 \iota_1 (\beta) = \begin{cases} \beta' & \mbox{ if $\beta$ has a bar that is not an anchor},\\
  \beta & \mbox{ if $\beta$ only has anchors.} \end{cases}
\]
From this definition we see the \textbf{exceptions} are those anchored barred permutations that have $k$ bars, all of which are anchors. (The exceptions form a subset of $B_{n,k}^0$.) Thus, from Observation \ref{obs3} the D.I.E. method has shown
\[
 \sum_{i=0}^k (-1)^i\binom{n+1}{i}(k+1-i)^n = \sum_{\beta \in \mathcal{B}_{n,k}} (-1)^{\# \textrm{float bars of } \beta} = \el{n}{k},
\]
which proves Equation \eqref{eq1}.

Before we move on, let us make some final observations. First, notice that the involution $\iota_1$ preserves the underlying permutation. With this in mind, define $\mathcal{B}_{\pi,k}$ to be the set of anchored barred permutations with underlying permutation $\pi$. Then we can partition $\mathcal{B}_{n,k}$ as
\[
 \mathcal{B}_{n,k} = \bigcup_{\pi \in S_n} \mathcal{B}_{\pi,k},
\]
and refine our formula permutation-by-permutation if we like. We find the structure of $\mathcal{B}_{\pi,k}$, and hence the action of $\iota_1$, is very simple at this level.

\begin{observation}\label{obs:perm}
Let $k\geq 0$ and let $\pi \in S_n$. If $\des(\pi)>k$, then by Observation \ref{obs2}, $\mathcal{B}_{\pi,k} =\emptyset$. If $\des(\pi) = k$, $\mathcal{B}_{\pi,k}$ contains one element with $k$ anchors and no other bars. If $\des(\pi)<k$, then every anchored barred permutation in $\mathcal{B}_{\pi,k}$ has at least one bar that is not an anchor, and hence no exceptional elements. Therefore,
\[
 \sum_{\beta \in \mathcal{B}_{\pi,k}} (-1)^{\# \textrm{float bars of } \beta} = \begin{cases} 1 & \mbox{ if } \des(\pi) = k,\\
 0 & \mbox{ otherwise.}
 \end{cases}
\]
\end{observation} 

\subsection{An alternating sign formula for $P$-Eulerian numbers}

Up to now, we have placed no restrictions on how our labeled balls can be placed into boxes. But imagine, say, we want to declare that ball $3$ always goes to the left of ball $1$, or that ball $1$ is in the last nonempty box. Our generalization of Equation \eqref{eq1} accounts for these sorts of restrictions, which can be encoded in a partially ordered set, known as a \emph{poset} for short. This will lead us to the notion of a \emph{$P$-partition} and the $P$-Eulerian numbers. Once we have the definitions sorted out, we will get an alternating sign formula for $P$-Eulerian numbers directly from Observation \ref{obs:perm}.

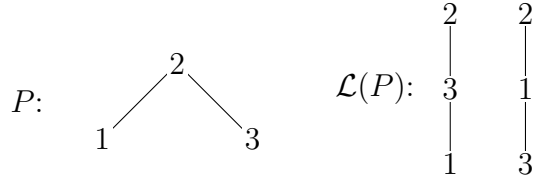
\begin{figure}
\[
 \begin{tikzpicture}[baseline=.5cm]
  \draw (0,1.5) node {$P$:};
  \draw (1,1) node[fill=white,inner sep=1] {$1$}-- (2,2) node[fill=white,inner sep=1] {$2$} -- (3,1) node[fill=white,inner sep=1] {$3$};
 \end{tikzpicture}
 \qquad
 \begin{tikzpicture}
  \draw (0,1) node {$\mathcal{L}(P)$:};
  \draw (1,0) node[fill=white,inner sep=1] {$1$}-- (1,1) node[fill=white,inner sep=1] {$3$} -- (1,2) node[fill=white,inner sep=1] {$2$};
  \draw (2,0) node[fill=white,inner sep=1] {$3$}-- (2,1) node[fill=white,inner sep=1] {$1$} -- (2,2) node[fill=white,inner sep=1] {$2$};
 \end{tikzpicture}
\]
\caption{The Hasse diagram of a poset $P$ and its linear extensions.}\label{fig:LP}
\end{figure}

Suppose we have $3$ balls and we want to declare that ball 2 is always the largest ball in the last nonempty box, i.e., that the corresponding barred permutation has 2 as the final number. We can capture this property by declaring that the balls are partially ordered as shown in Figure \ref{fig:LP}. This example shows $1<_P 2$ and $3<_P 2$. In general we illustrate such orderings with a graph known as a \emph{Hasse diagram}. The Hasse diagram has an edge from $i$ to $j$ for every pair $i<_P j$ such that if $i\leq_P x \leq_P j$, then $x = i$ or $x= j$. (All other comparisons follow by transitivity.) In pictures the edges are typically directed from bottom to top.

With a particular poset $P$, a pair of elements might or might not be comparable. In the example of Figure \ref{fig:LP}, elements 1 and 3 are incomparable. A poset with no incomparable elements is called a \emph{chain}, whereas a poset in which every pair of elements is incomparable is called an \emph{antichain}. The set of \emph{linear extensions} of a poset $P$, denoted $\LL(P)$, is the set of all chains that refine the partial ordering in the sense that if $i<_P j$, then $i$ is below $j$ in the chain. 

If $P$ is a partial ordering of $\{1,2,\ldots,n\}$, we readily identify its linear extensions with permutations. We write $i <_{\pi} j$ to mean $\pi^{-1}(i) < \pi^{-1}(j)$, i.e., $i$ appears to the left of $j$ in the one-line notation for $\pi$:
\[
 \pi(1) <_{\pi}  \pi(2) <_{\pi} \cdots <_{\pi} \pi(n).
\] 
We define the \emph{$P$-Eulerian numbers} $\el{P}{k}$ to be the number of linear extensions of $P$ with $k$ descents:
\[
 \el{P}{k} = |\{ \pi \in \LL(P) : \des(\pi) = k\}|.
\]

For the example in Figure \ref{fig:LP}, we have $\LL(P) = \{132, 312\}$. This gives $\el{P}{0} = 0, \el{P}{1}=2, \el{P}{2} = 0$. As a different example, if $P$ is an antichain, then there are no restrictions and every permutation is a linear extension. In this case, the $P$-Eulerian numbers are the usual Eulerian numbers: $\el{P}{k} = \el{n}{k}$.

Our next task is to explain how exactly to put balls into boxes with the restrictions from a poset $P$. We will see this amounts to considering only those barred permutations whose underlying permutation is a linear extension of $P$. Ultimately this gives an alternating sign formula for $P$-Eulerian numbers that generalizes Equation \eqref{eq1}. 

Throughout the rest of this section, fix a partial ordering $P$ on $\{1,2,\ldots,n\}$. A \emph{$P$-partition} is an order preserving function $f: \{1,2,\ldots,n\} \to \NN=\{0,1,2,\ldots\}$ such that for $i<_P j$:
\begin{itemize}
 \item $f(i) \leq f(j)$ if $i< j$,
 \item $f(i) < f(j)$ if $i > j$.
\end{itemize}
The set of $P$-partitions is denoted $\mathcal{A}(P)$.

Returning to the example of Figure \ref{fig:LP}, we see every $P$-partition $f$ must satisfy
\[
 f(1)\leq f(2) > f(3),
\]
so we have
\[
 \mathcal{A}(P) = \{\, (a_1, a_2, a_3)\in  \NN^3 : a_1 \leq a_2 > a_3 \,\}
\]
where we are identifying each $P$-partition $f$ with the values $(f(1),f(2),f(3))$. We can write $\A(P)$ as the disjoint union
\[
 \A(P) = \{ a_1 \leq a_3 < a_2 \} \cup \{ a_3 < a_1 \leq a_2 \},
\]
where each of these smaller sets can viewed as the set of $P$-partitions for a chain:
\[
 \A(P) = \A(132) \cup \A(312).
\]

By induction on the number of incomparable pairs in a poset $P$, we can see that the set of all $P$-partitions is the disjoint union of the sets of $P$-partitions for its linear extensions:
\begin{equation}\label{eq:Ap}
 \A(P) = \bigcup_{\pi \in \LL(P)} \A(\pi).
\end{equation}

The connection between $P$-partitions and placing labeled balls in labeled boxes is straightforward. A $P$-partition is nothing but a list of instructions for where to place the balls, with $f(i)$ recording the label of the box in which ball $i$ is placed. In this way, each $P$-partition can be identified with a barred permutation. For example, the $P$-partition
\[
(f(1),f(2),f(3),f(4),f(5),f(6),f(7),f(8),f(9))=(5,5,2,6,3,3,7,7,7)
\]
corresponds to $\beta=||3|56||12|4|789$. We say a barred permutation is \emph{$P$-compatible} if it arises in this way.

Define $\op_P(k)$ to be the number of $P$-partitions such that $f(i)\leq k$ for all $i\in P$. Equivalently, $\op_P(k)$ is the number of $P$-compatible barred permutations with $k$ bars. Notice that if $P$ is the antichain on $\{1,2,\ldots,n\}$ then $\op_P(k) = (k+1)^n$ as in Observation \ref{obs1}. 

Our generalization of Equation \eqref{eq1} is the following:
\begin{equation}\label{eq:Peul}
 \el{P}{k} = \sum_{i=0}^k (-1)^i\binom{n+1}{i}\op_P(k-i).
\end{equation}

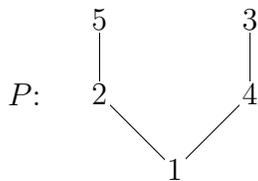
\begin{figure}
\[
 \begin{tikzpicture}[baseline=.5cm]
  \draw (0,2) node {$P$:};
  \draw (1,3) node[fill=white,inner sep=1] {$5$} -- (1,2) node[fill=white,inner sep=1] {$2$} -- (2,1) node[fill=white,inner sep=1] {$1$} -- (3,2) node[fill=white,inner sep=1] {$4$} -- (3,3) node[fill=white,inner sep=1] {$3$};
 \end{tikzpicture}
\]
\caption{The Hasse diagram of a poset $P$.}\label{fig:example}
\end{figure}

As an example, take $n = 5$ and the poset $P$ shown in Figure~\ref{fig:example}. There are six linear extensions: $12435$, $12453$, $12543$, $14235$, $14253$, and $14325$, with the number of descents: $1$, $1$, $2$, $1$, $2$, and $2$, respectively. A $P$-partition must satisfy $0 \leq f(1) \leq f(2) \leq f(5)$ and $0 \leq f(1) \leq f(4) < f(3)$. If $f(i) \leq k$ for all $i \in P$ and $f(1) = j$, we have $\binom{k-j+2}2$ choices for $(f(2),f(5))$ and $\binom{k-j+1}2$ choices for $(f(4),f(3))$. That means that
$$\op_P(k) = \sum_{j=0}^{k-1} \binom{k-j+2}2 \binom{k-j+1}2 = \frac{k(k+1)(k+2)(k+3)(2k+3)}{40}.$$
Indeed,
$$\sum_{i=0}^k (-1)^i\binom{6}{i}\op_P(k-i) = \begin{cases} 3, & k \in \{1,2\}\\
            0, & \text{otherwise}
		 \end{cases}.$$

We mirror our D.I.E. proof for Equation \eqref{eq1} and let $B_{P,k}^i$ denote the set of $P$-compatible anchored barred permutations with $k$ bars, $i$ of which floats. Next define $\mathcal{B}_{P,k}$ to be the set of all anchored barred permutations compatible with $P$ that have $k$ bars. We have both
\[
 \mathcal{B}_{P,k} = \bigcup_{i=0}^k B_{P,k}^i,
\]
and in light of \eqref{eq:Ap}, we also have
\begin{equation}\label{eq:Bp}
 \mathcal{B}_{P,k} = \bigcup_{\pi \in \LL(P)} \mathcal{B}_{\pi,k}.
\end{equation}

Mimicking Observation \ref{obs4} we find $|B_{P,k}^i| = \binom{n+1}{i}\op_P(k-i)$, and so our \textbf{description} for Equation \eqref{eq:Peul} is the set $\mathcal{B}_{P,k}$. But by Observation \ref{obs:perm} we and Equation \eqref{eq:Bp} we can apply the \textbf{involution} $\iota_1$ to find:
\begin{align*}
 \sum_{i=0}^k (-1)^i\binom{n+1}{i}\op_P(k-i) &= \sum_{\beta \in \mathcal{B}_{P,k}} (-1)^{\# \textrm{float bars of } \beta},\\
 &=\sum_{\pi \in \LL(P)} \sum_{\beta \in \mathcal{B}_{\pi,k}} (-1)^{\# \textrm{float bars of } \beta},\\
 &=\el{P}{k},
\end{align*}
which proves Equation \eqref{eq:Peul}.

We encourage the interested reader to learn more about $P$-partitions by visiting \cite[Section 3.15]{Stanley} or \cite[Section 3.2]{Petersen} and the references therein. Interestingly, $P$-partitions can be interpreted as the integer points in dilations of \emph{order polytopes} and $\op_P(k)$ counts the number of integer points. For example, when $P$ is an antichain of $n$ points, the order polytope is the hypercube $[0,1]^n$ and $(k+1)^n$ is the number of integer points in its $k$-fold dilation. This is a small corner of the broader research topic known as ``Ehrhart Theory'' which studies topics like the number of lattice points in rational polytopes.  

\section{Combinatorial topology}

We now focus our attention on the alternating sign formula in Equation \eqref{eq2} and its D.I.E. proof. The generalization of this approach lends itself to the study of abstract simplicial complexes.

\subsection{Second D.I.E. proof}

For this proof we use the model of set compositions. A \emph{set composition} is an ordered set partition of the set $\{1,2,\ldots,n\}$. For example, here is a composition of $\{1,2,\ldots,9\}$ with five blocks:
\[
 ( \{3\}, \{5,6\}, \{1,2\}, \{4\}, \{7,8,9\} )
\]
which we abbreviate with notation similar to barred permutations:
\[
 \gamma = 3|56\anc 12|4|789.
\]
In fact, we can think of set compositions as barred permutations with no floats and at most one bar per gap, i.e., no empty boxes. Importantly, we still mark each descent with an anchor.

For each partition of an $n$-element set into $k$ blocks, we can choose to order the blocks in $k!$ ways. Therefore, as $\st{n}{k}$ counts the number of set partitions of an $n$-element set into $k$ parts, the number of set compositions into $k$ parts is given by $\st{n}{k}k!$. These are sometimes known as the \emph{ordered Stirling numbers}.

A \emph{decorated set composition} is one in which some number of elements are highlighted, and we require that the maximal element in each block is highlighted. We indicate this with a horizontal bar over the elements, e.g.,
\[
\bar 3|5\bar 6\anc \bar 1\bar 2| \bar 4| \bar 78\bar 9
\]
is a decorated set composition.  The example above has five blocks, and seven highlighted elements. Two of these are non-maximal elements in their blocks: $1$ and $7$. We call highlighted elements that are non-maximal in this way \emph{remarkable elements} of $\gamma$. The sign of a decorated composition $\gamma$ is $(-1)^{\#\textrm{ remarkable elements }}$.

Let $C_{n,k}^i$ denote the set of decorated set compositions with $k+1$ highlighted elements and $i$ blocks. Notice that a decorated set composition of $n$ with $k+1$ highlighted elements can have at most $k+1$ blocks. Thus if we fix $k$, then $1\leq i \leq k+1$. If $\gamma$ is an element of $C_{n,k}^i$, each of the $i$ maximal elements is highlighted. This leaves $n-i$ elements from which to choose the remaining $k+1-i$ highlighted (remarkable) elements. Thus we make the following observation.

\begin{observation}\label{obs5}
Let $1\leq i \leq k+1 \leq n$. Then $|C_{n,k}^i| = \binom{n-i}{k+1-i}\st{n}{i}i!$. Moreover, each $\gamma$ in $C_{n,k}^i$ has sign $(-1)^{k+1-i}$.
\end{observation}

Therefore, our \textbf{description} for the sum in Equation \eqref{eq2} is given by the set $\mathcal{C}_{n,k} = \bigcup_{i=1}^{k+1} C_{n,k}^i$.

To define our involution, we first identify the leftmost highlighted element of $\gamma$ that does not have an anchor immediately to its right. Call this element the \emph{toggle element} of $\gamma$. If it has a vertical bar to its right we remove the bar. If it has no vertical bar to its right, we add that bar. (We exclude the maximal element in the rightmost block from consideration as a toggle, since it can never have a vertical bar to its right.) The new decorated set composition is $\gamma'$. Here are three examples with $n=9$ and $k=5$:
\[
\begin{array}{rcl}
\\
 3\bar 5 \bar 6\anc \bar 12\bar4|7\bar 8|\bar 9 & \leftrightarrow & 3\bar 5| \bar 6\anc \bar 12\bar4|7\bar 8|\bar 9 \\
\\
 3 5 \bar 6\anc \bar 1 |\bar 2 \bar 4\bar 7|8\bar 9 & \leftrightarrow & 3 5 \bar 6\anc \bar 1 \bar 2 \bar 4\bar 7|8\bar 9 \\
 \\
  5 \bar 9\anc \bar 3 \anc  2 \bar 4\anc \bar 1\bar 8 \anc 6 \bar 7 & \leftrightarrow &   5 \bar 9\anc \bar 3 \anc  2 \bar 4\anc \bar 1 | \bar 8 \anc 6 \bar 7 \\
 {}\\
\end{array}
\]
Importantly, we note that $\gamma$ and $\gamma'$ have the same toggle element. In one of them that element is largest in its block, while in the other it is a remarkable element. Thus the number of remarkable elements of $\gamma'$ is exactly one more or less than the number of remarkable elements of $\gamma$.

The \textbf{involution} $\iota_2: \mathcal{C}_{n,k} \to \mathcal{C}_{n,k}$ is as follows:
\[
 \iota_2 (\gamma) = \begin{cases} \gamma' & \mbox{ if $\gamma$ has a toggle element},\\
  \gamma & \mbox{ if $\gamma$ has no toggle element.} \end{cases}
\]
From this definition, we see the \textbf{exceptions} are those decorated set compositions for which there is no toggle. This means that every highlighted element (other than the rightmost element) is immediately followed by an anchor and there are no remarkable elements. Elements that are followed by anchors correspond to descents of the underlying permutation. So if $k+1$ is the number of highlighted elements, the underlying permutation has $k$ descents. Thus, the D.I.E. method has shown
\[
 \sum_{i=1}^{k+1} (-1)^{k+1-i}\binom{n-i}{k+1-i}\st{n}{i}i! = \sum_{\gamma \in \mathcal{C}_{n,k}} (-1)^{\# \textrm{remarkable elements of } \gamma} = \el{n}{k},
\]
as desired.

We now proceed to our topological generalization of this result.

\subsection{Combinatorial topology background}

We can view the result of the previous section from the more general perspective of abstract simplicial complexes. An \emph{abstract simplicial complex}, denoted $\Sigma$, is a collection  of sets called \emph{faces}, such that every subset of a face is a face. That is, if $F \in \Sigma$ and $G \subset F$, then $G \in \Sigma$. The \emph{dimension} of a face is one less than the cardinality of the face: $\dim(F) = |F|-1$. The zero dimensional faces are called \emph{vertices}, and every face is expressible as a union of vertices. The dimension of $\Sigma$ itself is the maximum dimension of a face in $\Sigma$: 
\[
\dim(\Sigma) = \max\{ \dim(F) : F \in \Sigma\}.
\]
Every abstract simplicial can be realized geometrically via points for vertices, edges for one-dimensional faces, triangles for two-dimensional faces, tetrahedra for three-dimensional faces, and so on. We will not bother with the details of how to embed a simplicial complex into Euclidean space, but for small examples we can draw pictures and use our visual understanding.

\begin{figure}
\[
\begin{array}{cc}
 \begin{tikzpicture}[scale=1.5]
  \draw[line width=1, fill=gray!40!white] (0,0) node[fill=black,circle,inner sep=1]{} --(1,1) node[fill=black,circle,inner sep=1]{} --(1,0) node[fill=black,circle,inner sep=1]{} --(0,0);
  \draw[line width=1] (1,1)--(2,1) node[fill=black,circle,inner sep=1]{} --(2,0) node[fill=black,circle,inner sep=1]{} --(1,0)--(1,1);
  \draw (0,0) node[left] {$a$};
  \draw (1,0) node[below] {$b$};
  \draw (1,1) node[above left] {$c$};
  \draw (2,1) node[above right] {$e$};
  \draw (2,0) node[right] {$d$};
 \end{tikzpicture}
&
\qquad
 \begin{tikzpicture}[xscale=1.1]
  \coordinate (0) at (0,0);
  \coordinate (a) at (-2,1);
  \coordinate (b) at (-1,1);
  \coordinate (c) at (0,1);
  \coordinate (d) at (2,1);
  \coordinate (e) at (1,1);
  \coordinate (ab) at (-2,2);
  \coordinate (ac) at (-1,2);
  \coordinate (bc) at (0,2);
  \coordinate (be) at (1,2);
  \coordinate (cd) at (2,2);
  \coordinate (de) at (3,2);
  \coordinate (abc) at (-1,3);
  \draw (0)--(a)--(ab)--(abc);
  \draw (a)--(ac);
  \draw (0)--(b)--(ab);
  \draw (0)--(c)--(ac)--(abc);
  \draw (0)--(d)--(cd);
  \draw (0)--(e)--(be);
  \draw (b)--(bc)--(abc);
  \draw (b)--(be);
  \draw (c)--(bc);
  \draw (c)--(cd);
  \draw (d)--(cd);
  \draw (d)--(de);
  \draw (e)--(be);
  \draw (e)--(de);
  \draw (0) node[fill=white, inner sep=1] {$\emptyset$};
  \draw (a) node [fill=white, inner sep=1] {$\{a\}$};
  \draw (b) node[fill=white, inner sep=1]  {$\{b\}$};
  \draw (c) node[fill=white, inner sep=1] {$\{c\}$};
  \draw (e) node[fill=white, inner sep=1] {$\{d\}$};
  \draw (d) node[fill=white, inner sep=1] {$\{e\}$};
  \draw (ab) node[fill=white, inner sep=1] {$\{a,b\}$};
  \draw (ac) node[fill=white, inner sep=1] {$\{a,c\}$};
  \draw (bc) node[fill=white, inner sep=1] {$\{b,c\}$};
  \draw (be) node[fill=white, inner sep=1] {$\{b,d\}$};
  \draw (cd) node[fill=white, inner sep=1] {$\{c,e\}$};
  \draw (de) node[fill=white, inner sep=1] {$\{d,e\}$};
  \draw (abc) node[fill=white, inner sep=1] {$\{a,b,c\}$};
 \end{tikzpicture} 
 \\
\textrm{(a)} & \textrm{(b)}
 \end{array}
\]
\caption{Two visual representations of a nonpure simplicial complex.}\label{fig:nonpure}
\end{figure}
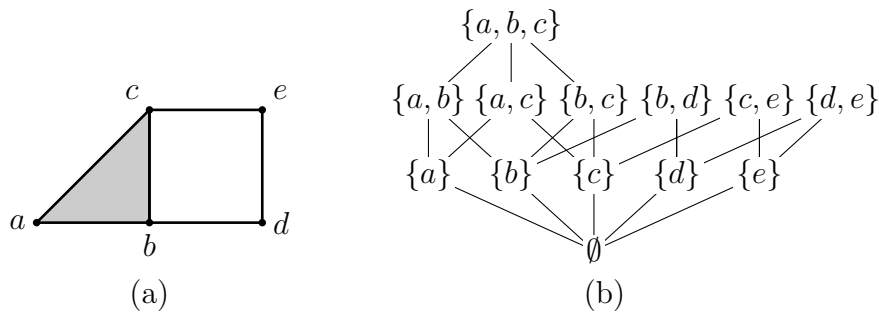

For example, Figure \ref{fig:nonpure}(a) shows a picture of a complex with five vertices, six edges, and one triangle. We can also visualize a simplicial complex $\Sigma$ as a partially ordered set with faces ordered by containment, as shown if Figure \ref{fig:nonpure}(b).

To any finite simplicial complex $\Sigma$ of dimension $d-1$, there are two well-studied invariants known as the \emph{$f$-vector}, denoted $f(\Sigma) = (f_0, f_1, f_2, \ldots, f_d)$ and the \emph{$h$-vector} $h(\Sigma) = (h_0, h_1, h_2, \ldots, h_d)$. The $f$-vector counts the faces of $\Sigma$ according to dimension, i.e., $f_i$ is the number of $(i-1)$-dimensional faces. The $h$-vector is a linear transformation of the $f$-vector given by
\begin{equation}\label{eq:ftoh}
 h_k = \sum_{i=0}^k (-1)^{k-i}\binom{d-i}{k-i} f_i.
\end{equation}
Note that $h_d = (-1)^d(f_0 - f_1 + f_2 -\cdots) = (-1)^d(1-\chi(\Sigma))$, where $\chi(\Sigma)$ is the Euler characteristic of $\Sigma$. For the example in Figure \ref{fig:nonpure} we can compute $f(\Sigma) = (1,5,6,1)$ and $h(\Sigma) = (1,2,-1,-1)$. 

One might wonder if there is a general D.I.E. explanation for Equation \eqref{eq:ftoh}, phrased in the language of simplicial complexes. As we will explain, the answer is an emphatic \textbf{yes} when we assume two properties of our complex. From now on, assume our complex $\Sigma$ is \emph{pure} and \emph{partitionable}.\footnote{There is also a commonly studied and slightly stronger property of ``shellability'' for a simplicial complex. See Chapter 8 of \cite{Ziegler} for an introduction to these ideas.}

To be pure simply means that all maximal faces have the same dimension. The example in Figure \ref{fig:nonpure} is not pure since it has a triangle as well as edges that are not contained in any triangle. In the case where $\Sigma$ is pure, we refer to the maximal faces as \emph{facets}.

To be \emph{partitionable} means that $\Sigma$ can be expressed as a disjoint union of intervals
\[
 [\widetilde{F},F] = \{ G : \widetilde{F} \subseteq G \subseteq F \},
\]
where each $F$ is a facet. Each such interval is a boolean lattice, isomorphic to the lattice of subsets $F-\widetilde{F}$. A partition of $\Sigma$ matches each facet $F$ with one of its lower-dimensional faces $\widetilde{F}$. Taking Equation \eqref{eq:ftoh} as the definition of the $h$-vector, we will show that 
\begin{equation}\label{eq:hk}
 h_k = |\{ \textrm{facets } F \in \Sigma : |\widetilde{F}| = k\}|.
\end{equation}
It is noteworthy that this result does not depend on the choice of partition, only the fact that a partition exists for $\Sigma$.

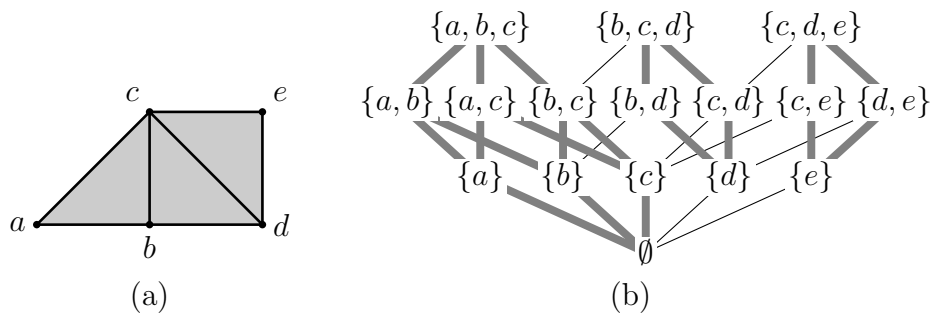
\begin{figure}
\[
\begin{array}{cc}
 \begin{tikzpicture}[scale=1.5]
  \draw[fill=gray!40!white, line width=1] (0,0) node[fill=black,circle,inner sep=1]{} --(1,1) node[fill=black,circle,inner sep=1]{} --(2,1) node[fill=black,circle,inner sep=1]{} --(2,0) node[fill=black,circle,inner sep=1]{} --(1,0) node[fill=black,circle,inner sep=1]{} --(0,0);
  \draw[line width=1] (1,1)--(2,0);
  \draw[line width=1] (1,1)--(1,0);
  \draw (0,0) node[left] {$a$};
  \draw (1,0) node[below] {$b$};
  \draw (1,1) node[above left] {$c$};
  \draw (2,1) node[above right] {$e$};
  \draw (2,0) node[right] {$d$};
 \end{tikzpicture}
&
\quad
 \begin{tikzpicture}[xscale=1.1]
  \coordinate (0) at (0,0);
  \coordinate (a) at (-2,1);
  \coordinate (b) at (-1,1);
  \coordinate (c) at (0,1);
  \coordinate (d) at (1,1);
  \coordinate (e) at (2,1);
  \coordinate (ab) at (-3,2);
  \coordinate (ac) at (-2,2);
  \coordinate (bc) at (-1,2);
  \coordinate (bd) at (0,2);
  \coordinate (cd) at (1,2);
  \coordinate (ce) at (2,2);
  \coordinate (de) at (3,2);
  \coordinate (abc) at (-2,3);
  \coordinate (bcd) at (0,3);
  \coordinate (cde) at (2,3);
  \draw (0)--(a)--(ab)--(abc);
  \draw (a)--(ac);
  \draw (0)--(b)--(ab);
  \draw (0)--(c)--(ac)--(abc);
  \draw (0)--(d)--(cd)--(cde);
  \draw (0)--(e)--(de);
  \draw (b)--(bc)--(abc);
  \draw (b)--(bd)--(bcd);
  \draw (c)--(bc)--(bcd);
  \draw (c)--(cd)--(bcd);
  \draw (c)--(ce);
  \draw (d)--(bd);
  \draw (d)--(de)--(cde);
  \draw (e)--(ce)--(cde);
  \draw[rounded corners, line width=3, gray] (0)--(a)--(ab)--(abc)--(ac)--(a);
  \draw[rounded corners, line width=3, gray] (0)--(b)--(ab);
  \draw[rounded corners, line width=3, gray] (0)--(c)--(bc)--(abc);
  \draw[rounded corners, line width=3, gray] (b)--(bc);
  \draw[rounded corners, line width=3, gray] (c)--(ac);
  \draw[rounded corners, line width=3, gray] (d)--(bd)--(bcd)--(cd)--(d);
  \draw[rounded corners, line width=3, gray] (e)--(ce)--(cde)--(de)--(e);
  \draw (0) node[fill=white, inner sep=1] {$\emptyset$};
  \draw (a) node [fill=white, inner sep=1] {$\{a\}$};
  \draw (b) node[fill=white, inner sep=1]  {$\{b\}$};
  \draw (c) node[fill=white, inner sep=1] {$\{c\}$};
  \draw (d) node[fill=white, inner sep=1] {$\{d\}$};
  \draw (e) node[fill=white, inner sep=1] {$\{e\}$};
  \draw (ab) node[fill=white, inner sep=1] {$\{a,b\}$};
  \draw (ac) node[fill=white, inner sep=1] {$\{a,c\}$};
  \draw (bc) node[fill=white, inner sep=1] {$\{b,c\}$};
  \draw (bd) node[fill=white, inner sep=1] {$\{b,d\}$};
  \draw (ce) node[fill=white, inner sep=1] {$\{c,e\}$};
  \draw (cd) node[fill=white, inner sep=1] {$\{c,d\}$};
  \draw (de) node[fill=white, inner sep=1] {$\{d,e\}$};
  \draw (abc) node[fill=white, inner sep=1] {$\{a,b,c\}$};
  \draw (bcd) node[fill=white, inner sep=1] {$\{b,c,d\}$};
  \draw (cde) node[fill=white, inner sep=1] {$\{c,d,e\}$};
 \end{tikzpicture} 
 \\
\textrm{(a)} & \textrm{(b)}
 \end{array}
\]
\caption{A pure partitionable simplicial complex.}\label{fig:pure}
\end{figure}

In Figure \ref{fig:pure} we see an example of pure partitionable complex $\Sigma$ of dimension two. We see that $f(\Sigma) = (1,5,7,3)$ and $h(\Sigma) = (1,2,0,0)$. From Figure \ref{fig:pure} (b) we show that $\Sigma$ can be partitioned into three intervals:
\[
 \Sigma = [\emptyset, \{a,b,c\}] \cup [\{d\}, \{b,c,d\}] \cup [\{e\},\{c,d,e\}].
\]
One of these has $|\widetilde{F}|=0$ and two of them have $|\widetilde{F}| = 1$, which affirms $h_0=1$ and $h_1 = 2$.

As we shall later see, the set of all set compositions of $\{1,2,\ldots,n\}$ can be identified with faces of a partitionable simplicial complex with $h_k = \el{n}{k}$, so Equation \eqref{eq:ftoh} truly is a generalization of Equation \eqref{eq2}.

\subsection{The topological D.I.E. argument}

We will now describe how to apply the D.I.E. method to this situation. Let $\Sigma$ be a finite, pure, partitionable complex. We will fix a partition of $\Sigma$, and for reasons that will be important later, we choose a labeling of the vertices: $v_1, \ldots, v_n$. This choice of labeling, like the choice of partition, is arbitrary, but must be fixed for the argument to make sense.

Once we have fixed our partition, every face $G \in \Sigma$ belongs to a unique interval $[\widetilde{F},F]$. We say $\widetilde{F} = \widetilde{F}(G)$ is the \emph{anchor set} of $G$. Notice that $G = \widetilde{F}\cup I$ for some subset $I \subseteq F-\widetilde{F}$. 

We now define a \emph{decorated face} of $\Sigma$ to be a pair $(G,J)$, denoted $G^J$ for short, where $G \in [\widetilde{F},F]$ and $J \subseteq F-G$. We refer to $J$ as the set of \emph{remarkable vertices} of the decorated face. The sign of a decorated face $G^J$ is $(-1)^J$.

Let $C_{\Sigma,k}^i$ denote the set of decorated faces such that $|G|=i$ and $|J| = k-i$. If $|G|=i$, then $|F-G|=d-i$, so there are $\binom{d-i}{k-i}$ ways to choose the remarkable elements $J \subseteq F-G$. Thus we can make the following observation that mimics Observation \ref{obs5}.

\begin{observation}\label{obs6}
Let $\Sigma$ be a pure partitionable simplicial complex of dimension $(d-1)$, and let $f=(f_0,f_1,\ldots,f_d)$ denote its $f$-vector. Then for $0\leq i \leq k \leq d$, $|C_{\Sigma,k}^i| = \binom{d-i}{k-i}f_i$. Moreover, each decorated face in $C_{\Sigma,k}^i$ has sign $(-1)^{k-i}$. 
\end{observation} 

Therefore the \textbf{description} for the sum in Equation \eqref{eq:ftoh} is given by the set $\mathcal{C}_{\Sigma,k} = \bigcup_{i=0}^k C_{\Sigma,k}^i$. In other words, $\mathcal{C}_{\Sigma,k}$ is the set of decorated faces $G^J$ with $|G \cup J| = k$.

To define our involution, we consider a decorated face $G^J$. We let $I = G-\widetilde{F}$ and let $v^* \in I \cup J$ denote the vertex of least index in $I \cup J$. (Recall that the vertices of $\Sigma$ are labeled, so this notion is well-defined.) Our toggle is to simply move $v^*$ betweeen sets $I$ and $J$. That is, if $v^* \in J$, let $J'=J-\{v^*\}$ and $G'=G\cup \{v^*\}$. If $v^* \notin J$, then $J'=J\cup \{v^*\}$ and $G'=G-\{v^*\}$. Notice that $G^J$ and $G'^{J'}$ are both members of $\mathcal{C}_{\Sigma,k}$ since $|G\cup J| = |G'\cup J'|=k$, and they have opposite sign since the cardinalities of $J$ and $J'$ differ by one. 

Thus the \textbf{involution} $\iota_3: \mathcal{C}_{\Sigma,k} \to \mathcal{C}_{\Sigma,k}$ is as follows:
\[
 \iota_3(G^J) = \begin{cases} 
  G'^{J'} & \mbox{ if $I\cup J$ is nonempty},\\
  G^J & \mbox{ if $I=J=\emptyset$.}
 \end{cases}
\]

If $I = \emptyset$, then $G=\widetilde{F}$, so from the definition of $\iota_3$, the exceptions are those decorated faces of the form $\widetilde{F}^{\emptyset}$ with $|\widetilde{F}| = k$, each of which has positive sign. Therefore the D.I.E. method shows that
\begin{equation}\label{eq:DIEsimp}
 h_k = \sum_{i=0}^k (-1)^{k-i}\binom{d-i}{k-i} f_i = \sum_{G^J \in \mathcal{C}_{\Sigma,k}} (-1)^J = |\{ \textrm{facets } F \in \Sigma : |\widetilde{F}| = k\}|,
\end{equation}
as desired.

\subsection{Back to Eulerian numbers}

We now describe how set compositions arise in a topological context and how Equation \eqref{eq:ftoh} gives rise to Equation \eqref{eq3} for Eulerian numbers. To do this, we will produce a simplicial complex $\Delta_n$ of dimension $(n-1)$ with $f$-vector $f=(f_0,f_1,\ldots,f_{n})$ such that  $f_i = \st{n+1}{i+1}i!$ for each $i=0,1,\ldots,n$. Moreover it has a partition such that the set of faces $\widetilde{F}$ with $|\widetilde{F}| = k$ is in bijection with the set of permutations in $S_n$ with $k$ descents. To make these connections we follow the presentation of Section 9.2 of \cite{Petersen}, which then ties back to Equation \eqref{eq2}.

To continue, we need a couple more definitions. First define a \emph{flag} of a simplicial complex $\Sigma$ to be a sequence of distinct faces $F_1, F_2, \ldots,F_k$ such that
\[
 F_1 \subset F_2 \subset \cdots \subset F_k.
\]
The \emph{barycentric subdivision} of an abstract simplicial complex $\Sigma$ is the complex whose faces are flags. In what follows we let $\Sigma_n$ be the simplex whose faces are all subsets of $\{1,2,\ldots,n\}$, and let $\Delta_n$ denote its barycentric subdivision.

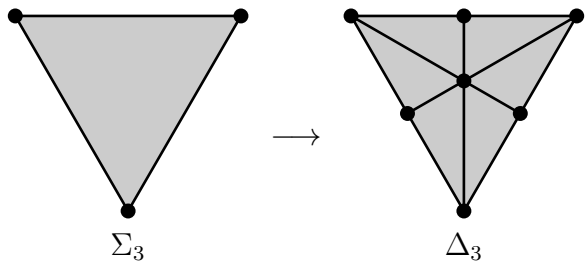
\begin{figure}
\[
\begin{array}{c}
\begin{tikzpicture}[scale=1.5]
\draw[fill=gray!40!white,line width=1] (0,0) -- (1,-1.73) -- (2,0) -- (0,0);
\draw (0,0) node[fill=black,circle,inner sep=2]{};
\draw (1,-1.73) node[fill=black,circle,inner sep=2]{};
\draw (2,0) node[fill=black,circle,inner sep=2]{};
\end{tikzpicture}
\\
\Sigma_3
\end{array}
\longrightarrow
\begin{array}{c}
\begin{tikzpicture}[scale=1.5]
\draw[fill=gray!40!white,line width=1] (0,0) -- (1,-1.73) -- (2,0) -- (0,0);
\draw[line width=1] (0,0) node[fill=black,circle,inner sep=2]{} -- (1.5,-1.73*0.5) node[fill=black,circle,inner sep=2]{};
\draw[line width=1] (1,0) node[fill=black,circle,inner sep=2]{} -- (1,-1.73) node[fill=black,circle,inner sep=2]{};
\draw[line width=1] (.5,-1.73*0.5) node[fill=black,circle,inner sep=2]{} -- (2,0) node[fill=black,circle,inner sep=2]{};
\draw (1,-1.73*0.333) node[fill=black,circle,inner sep=2]{};
\end{tikzpicture}
\\
\Delta_3
\end{array}
\]
\caption{A $2$-dimensional simplex $\Sigma_3$ and its barycentric subdivision $\Delta_3$.}\label{fig:barycenter}
\end{figure}

For example, listing the flags of $\Sigma_3$ we find the faces of $\Delta_3$:
\[
\begin{array}{c | c | c | c}
 \mbox{empty face} & \mbox{vertices} & \mbox{edges} & \mbox{triangles} \\
 \hline
\emptyset & \{1\} & \{1\} \subset \{1,2\} & \{1\} \subset \{1,2\} \subset \{1,2,3\} \\
 & \{ 2\} & \{1\} \subset \{1,3\} & \{1\} \subset \{1,3\} \subset \{1,2,3\} \\
 & \{3\} & \{2\} \subset \{1,2\} & \{2\} \subset \{1,2\} \subset \{1,2,3\} \\
 & \{1,2\} & \{2\} \subset \{2,3\} & \{2\} \subset \{2,3\} \subset \{1,2,3\} \\
 & \{1,3\} & \{3\} \subset \{1,3\} & \{3\} \subset \{1,3\} \subset \{1,2,3\} \\
 & \{2,3\} & \{3\} \subset \{2,3\} & \{3\} \subset \{2,3\} \subset \{1,2,3\} \\
 & \{1,2,3\} & \{1\} \subset \{1,2,3\} \\
 & & \{2\} \subset \{1,2,3\} \\
 & & \{3\} \subset \{1,2,3\} \\
 & & \{1,2\} \subset \{1,2,3\} \\
 & & \{1,3\} \subset \{1,2,3\} \\
 & & \{2,3\} \subset \{1,2,3\}
\end{array}
\]
Computing the $f$- and $h$-vectors we find $f(\Delta_3)=(1,7,12,6)$, $h(\Delta_3) = (1,4,1,0)$.
We can visualize this as in Figure \ref{fig:barycenter}.

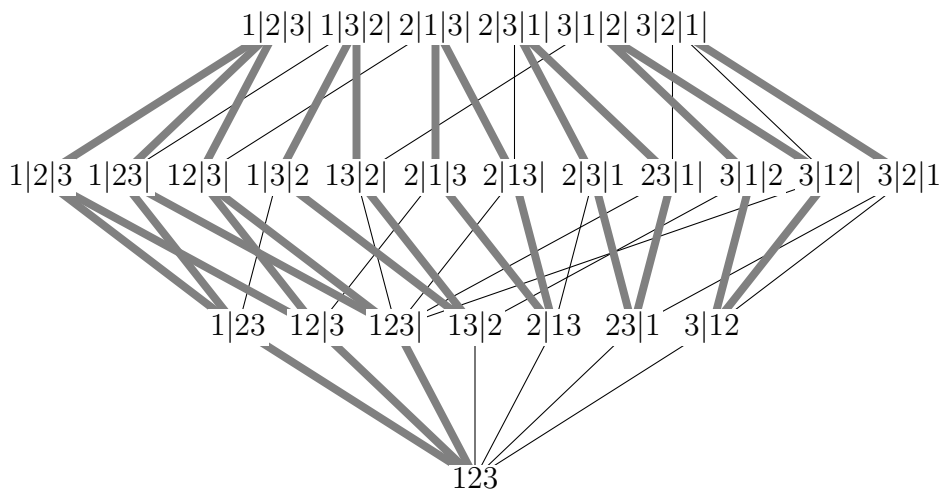
\begin{figure}
\[
 \begin{tikzpicture}[xscale=1.05,yscale=2]
  \coordinate (123) at (0,0);
  \coordinate (1b23) at (-3,1);
  \coordinate (12b3) at (-2,1);
  \coordinate (123b) at (-1,1);
  \coordinate (13b2) at (0,1);
  \coordinate (2b13) at (1,1);
  \coordinate (23b1) at (2,1);
  \coordinate (3b12) at (3,1);
  \coordinate (1b2b3) at (-5.5,2);
  \coordinate (1b23b) at (-4.5,2);
  \coordinate (12b3b) at (-3.5,2);
  \coordinate (1b3b2) at (-2.5,2);
  \coordinate (13b2b) at (-1.5,2);
  \coordinate (2b1b3) at (-0.5,2);
  \coordinate (2b13b) at (0.5,2);
  \coordinate (2b3b1) at (1.5,2);
  \coordinate (23b1b) at (2.5,2);
  \coordinate (3b1b2) at (3.5,2);
  \coordinate (3b12b) at (4.5,2);
  \coordinate (3b2b1) at (5.5,2);
  \coordinate (1b2b3b) at (-2.5,3);
  \coordinate (1b3b2b) at (-1.5,3);
  \coordinate (2b1b3b) at (-0.5,3);
  \coordinate (2b3b1b) at (0.5,3);
  \coordinate (3b1b2b) at (1.5,3);
  \coordinate (3b2b1b) at (2.5,3);
  \draw (123)--(1b23)--(1b2b3)--(1b2b3b);
  \draw (1b23)--(1b23b);
  \draw (123)--(12b3)--(1b2b3);
  \draw (123)--(123b)--(1b23b)--(1b2b3b);
  \draw (123)--(13b2)--(1b3b2)--(1b3b2b);
  \draw (123)--(2b13)--(2b1b3)--(2b1b3b);
  \draw (123)--(23b1)--(2b3b1)--(2b3b1b);
  \draw (123)--(3b12)--(3b1b2)--(3b1b2b);
  \draw (1b23)--(1b3b2);
  \draw (12b3)--(12b3b);
  \draw (12b3)--(2b1b3);
  \draw (123b)--(12b3b)--(1b2b3b);
  \draw (123b)--(2b13b);
  \draw (123b)--(3b12b);
  \draw (123b)--(23b1b);
  \draw (123b)--(13b2b);
  \draw (1b23b)--(1b3b2b);
  \draw (12b3b)--(2b1b3b);
  \draw (13b2)--(13b2b);
  \draw (13b2)--(3b1b2);
  \draw (13b2b)--(3b1b2b);
  \draw (2b13)--(2b13b)--(2b1b3b);
  \draw (2b13)--(2b3b1);
  \draw (2b13b)--(2b3b1b);
  \draw (23b1)--(23b1b)--(2b3b1b);
  \draw (23b1)--(3b2b1);
  \draw (23b1b)--(3b2b1b);
  \draw (3b12)--(3b12b)--(3b1b2b);
  \draw (3b12)--(3b2b1);
  \draw (3b12b)--(3b2b1b);
  \draw (3b2b1)--(3b2b1b);
  \draw (13b2b)--(1b3b2b);
  \draw[rounded corners, line width=3, gray] (123)--(1b23)--(1b2b3)--(1b2b3b);
  \draw[rounded corners, line width=3, gray] (1b23)--(1b23b);
  \draw[rounded corners, line width=3, gray] (123)--(12b3)--(1b2b3);
  \draw[rounded corners, line width=3, gray] (123)--(123b)--(1b23b)--(1b2b3b);
  \draw[rounded corners, line width=3, gray] (12b3)--(12b3b)--(1b2b3b);
  \draw[rounded corners, line width=3, gray] (123b)--(12b3b);
  \draw[rounded corners, line width=3, gray] (13b2)--(1b3b2)--(1b3b2b)--(13b2b)--(13b2);
  \draw[rounded corners, line width=3, gray] (2b13)--(2b1b3)--(2b1b3b)--(2b13b)--(2b13);
  \draw[rounded corners, line width=3, gray] (23b1)--(2b3b1)--(2b3b1b)--(23b1b)--(23b1);
  \draw[rounded corners, line width=3, gray] (3b12)--(3b1b2)--(3b1b2b)--(3b12b)--(3b12);
  \draw[rounded corners, line width=3, gray] (3b2b1)--(3b2b1b);
  \draw (123) node[fill=white, inner sep=1] {$123$};
  \draw (1b23) node [fill=white, inner sep=1] {$1|23$};
  \draw (12b3) node[fill=white, inner sep=1]  {$12|3$};
  \draw (123b) node[fill=white, inner sep=1] {$123|$};
  \draw (1b2b3) node[fill=white, inner sep=1] {$1|2|3$};
  \draw (1b23b) node[fill=white, inner sep=1] {$1|23|$};
  \draw (12b3b) node[fill=white, inner sep=1] {$12|3|$};
  \draw (1b2b3b) node[fill=white, inner sep=1] {$1|2|3|$};
  \draw (13b2) node[fill=white, inner sep=1] {$13|2$};
  \draw (1b3b2) node[fill=white, inner sep=1] {$1|3|2$};
  \draw (13b2b) node[fill=white, inner sep=1] {$13|2|$};
  \draw (1b3b2b) node[fill=white, inner sep=1] {$1|3|2|$};
  \draw (2b13) node[fill=white, inner sep=1] {$2|13$};
  \draw (2b1b3) node[fill=white, inner sep=1] {$2|1|3$};
  \draw (2b13b) node[fill=white, inner sep=1] {$2|13|$};
  \draw (2b1b3b) node[fill=white, inner sep=1] {$2|1|3|$};
  \draw (23b1) node[fill=white, inner sep=1] {$23|1$};
  \draw (2b3b1) node[fill=white, inner sep=1] {$2|3|1$};
  \draw (23b1b) node[fill=white, inner sep=1] {$23|1|$};
  \draw (2b3b1b) node[fill=white, inner sep=1] {$2|3|1|$};
  \draw (3b12) node[fill=white, inner sep=1] {$3|12$};
  \draw (3b1b2) node[fill=white, inner sep=1] {$3|1|2$};
  \draw (3b12b) node[fill=white, inner sep=1] {$3|12|$};
  \draw (3b1b2b) node[fill=white, inner sep=1] {$3|1|2|$};
  \draw (3b2b1) node[fill=white, inner sep=1] {$3|2|1$};
  \draw (3b2b1b) node[fill=white, inner sep=1] {$3|2|1|$};
 \end{tikzpicture} 
 \]
\caption{The face poset of $\Delta_3$, partitioned.}\label{fig:barycenterfaces}
\end{figure}

Now for the complex $\Delta_n$, we can simplify the notation for flags. Since flags are sequences of nested subsets of $\{1,2,\ldots,n\}$, we can choose to only keep track of when new elements are added. That is, given a flag of $k$ sets,
\[
  \emptyset \subseteq F_1 \subset F_2 \subset \cdots \subset F_k \subseteq \{1,2,\ldots,n\},
\]
we can write $A_i = F_i - F_{i-1}$, with $F_0 = \emptyset$ and $F_{k+1} = \{1,2,\ldots,n\}$. Instead of keeping track of each $F_i$, we can record the sequence of differences $A_1, A_2, \ldots, A_{k+1}$. For example, with $n=6$, the following is a flag with $k=3$ sets:
\[
  \emptyset \subseteq \{3\} \subset \{1,3,4\} \subset \{1,3,4,6\} \subseteq \{1,2,3,4,5,6\},
\]
and it has
\[
 (A_1,A_2,A_3,A_4) = (\{3\}, \{1,4\}, \{6\}, \{2,5\}).
\]
If we condense notation further, we can write this as a set composition: $3|14|6|25$. With this in mind, we can draw the face poset of $\Delta_3$ as shown in Figure \ref{fig:barycenterfaces}. Notice that if $\{1,2,3\}$ appears in a flag, then the final block of the composition is empty. Let $\overline{\mathcal{C}}_n$ denote the set of compositions of $\{1,2,\ldots,n\}$ such that the final block can be empty.

Let us count the faces by dimension. First, observe that $f_i$ is counting these set compositions according to the number of bars. We can see that any such composition with $i$ bars either has its final block empty or not. There are $\st{n}{i+1}(i+1)!$ compositions of $\{1,2,\ldots,n\}$ into $i+1$ nonempty blocks. If the final block is empty, then the first $i$ blocks form a composition of $\{1,2,\ldots,n\}$ into $i$ nonempty blocks, and there are $\st{n}{i}i!$ of these. In total, we have
\begin{align}
 f_i(\Delta_n) &= \st{n}{i+1}(i+1)! + \st{n}{i}i!, \nonumber \\
  &= i!\left(\st{n}{i+1}(i+1)+\st{n}{i}\right),  \nonumber \\
  &= i!\st{n+1}{i+1}, \label{eq:fiDelta}
\end{align}
where the final equality follows from a standard recurrence for the Stirling numbers.\footnote{To create a set partition of $n+1$ elements into $i+1$ blocks, we either add $n+1$ to a set partition of $n$ elements into $i+1$ blocks (in $i+1$ ways), or we add a new singleton block containing $n+1$ to a set partition of $n$ elements into $i$ blocks.} This shows that $\Delta_n$ has the $f$-vector we desire to make Equation \eqref{eq:DIEsimp} match Equation \eqref{eq3}.

Next we argue that $\Delta_n$ is partitionable and that counting its anchor sets by dimension is the same as counting permutations according to descents. But this is remarkably simple! We merely partition the faces of $\Delta_n$ according to the underlying permutation in the set composition, as indicated in Figure \ref{fig:barycenterfaces}. That is, each block $[\widetilde{F},F]$ corresponds to a permutation $\pi \in S_n$. It has as its maximal element the maximally refined set composition $F=\pi(1)|\pi(2)|\cdots|\pi(n)|$, and its minimal element $\widetilde{F}$ is the set composition in which bars only occur in descent positions of $\pi$. This allows us to make the following observation.

\begin{observation}
 The simplicial complex $\Delta_n$ is partitionable, with $|\{ \textrm{facets } F \in \Delta_n : |\widetilde{F}|=k\}| = \el{n}{k}$.
\end{observation}
 
Taking this observation together with Equation \eqref{eq:fiDelta} for $f_i(\Delta_n)$ and substituting in Equation \eqref{eq:DIEsimp} yields Equation \eqref{eq3}, our third and final alternating sum expression mentioned in the introduction.

\begin{figure}
\[
\begin{array}{cc}
 \begin{tikzpicture}[scale=1.5]
  \draw[line width=1] (0,0) node[fill=black,circle,inner sep=2]{} --(1,-1.73) node[fill=black,circle,inner sep=2]{} --(2,0) node[fill=black,circle,inner sep=2]{} --(1,0) node[fill=black,circle,inner sep=2]{}--(0,0);
  \draw (.5,-1.73*0.5) node[fill=black,circle,inner sep=2]{};
  \draw (1.5,-1.73*0.5) node[fill=black,circle,inner sep=2]{};
 \end{tikzpicture}
&
\qquad
 \begin{tikzpicture}[xscale=1.1,yscale=1.5]
  \coordinate (123) at (0,0);
  \coordinate (1b23) at (-2.5,1);
  \coordinate (12b3) at (-1.5,1);
  \coordinate (13b2) at (-0.5,1);
  \coordinate (2b13) at (0.5,1);
  \coordinate (23b1) at (1.5,1);
  \coordinate (3b12) at (2.5,1);
  \coordinate (1b2b3) at (-2.5,2);
  \coordinate (1b3b2) at (-1.5,2);
  \coordinate (2b1b3) at (-0.5,2);
  \coordinate (2b3b1) at (0.5,2);
  \coordinate (3b1b2) at (1.5,2);
  \coordinate (3b2b1) at (2.5,2);
  \draw (123)--(1b23)--(1b2b3);
  \draw (123)--(12b3)--(1b2b3);
  \draw (123)--(13b2)--(1b3b2);
  \draw (123)--(2b13)--(2b1b3);
  \draw (123)--(23b1)--(2b3b1);
  \draw (123)--(3b12)--(3b1b2);
  \draw (1b23)--(1b3b2);
  \draw (12b3)--(2b1b3);
  \draw (13b2)--(3b1b2);
  \draw (2b13)--(2b3b1);
  \draw (23b1)--(3b2b1);
  \draw (3b12)--(3b2b1);
  \draw[rounded corners, line width=3, gray] (123)--(1b23)--(1b2b3);
  \draw[rounded corners, line width=3, gray] (123)--(12b3)--(1b2b3);
  \draw[rounded corners, line width=3, gray] (13b2)--(1b3b2);
  \draw[rounded corners, line width=3, gray] (2b13)--(2b1b3);
  \draw[rounded corners, line width=3, gray] (23b1)--(2b3b1);
  \draw[rounded corners, line width=3, gray] (3b12)--(3b1b2);
  \draw (123) node[fill=white, inner sep=1] {$123$};
  \draw (1b23) node [fill=white, inner sep=1] {$1|23$};
  \draw (12b3) node[fill=white, inner sep=1]  {$12|3$};
  \draw (1b2b3) node[fill=white, inner sep=1] {$1|2|3$};
  \draw (13b2) node[fill=white, inner sep=1] {$13|2$};
  \draw (1b3b2) node[fill=white, inner sep=1] {$1|3|2$};
  \draw (2b13) node[fill=white, inner sep=1] {$2|13$};
  \draw (2b1b3) node[fill=white, inner sep=1] {$2|1|3$};
  \draw (23b1) node[fill=white, inner sep=1] {$23|1$};
  \draw (2b3b1) node[fill=white, inner sep=1] {$2|3|1$};
  \draw (3b12) node[fill=white, inner sep=1] {$3|12$};
  \draw (3b1b2) node[fill=white, inner sep=1] {$3|1|2$};
  \draw (3b2b1) node[fill=white, inner sep=1] {$3|2|1$};
 \end{tikzpicture} 
 \\
\textrm{(a)} & \textrm{(b)}
 \end{array}
\]
\caption{The complex $\Delta'_3$ and its face poset.}\label{fig:Delta'}
\end{figure}
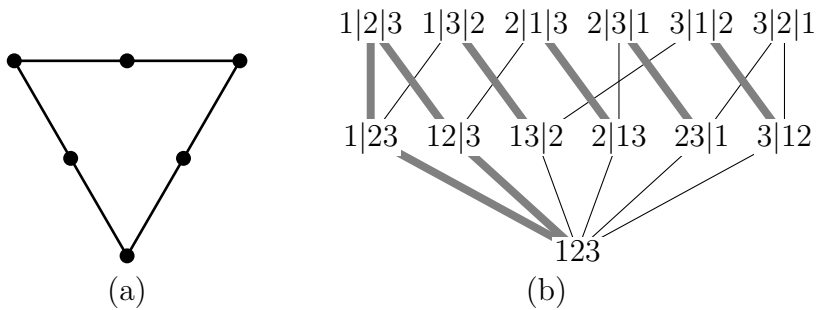

We conclude by pointing out that Equation \eqref{eq2} can also be argued via \eqref{eq:DIEsimp}. To this end, we consider the \emph{boundary complex} $\partial\Sigma_n$ of the simplex $\Sigma_n$. This is the simplicial complex whose faces are all \emph{proper} subsets $F \subsetneq \{1,2,\ldots,n\}$. We let $\Delta'_n$ denote the barycentric subdivision of $\partial\Sigma_n$. As we consider faces of $\Delta_n'$, i.e., flags in $\partial\Sigma_n$, we use the same bookkeeping strategy as $\Sigma_n$. The only difference now is that we do not allow $F_k = \{1,2,\ldots,n\}$ in a flag, and therefore the corresponding set compositions have no empty blocks. The complex $\Delta_3'$ and its face poset are shown in Figure \ref{fig:Delta'}.

Now the number of $i$-dimensional faces is the number of compositions with $i+1$ blocks, so $f_i(\Delta'_n) = \st{n}{i+1}(i+1)!$ and, just as with $\Delta_n$, we can partition the faces of $\Delta_n'$ into blocks according to the underlying permutation, yielding $h_k = \el{n}{k}$ again. Putting these facts together with Equation \eqref{eq:ftoh} yields the identity of Equation \eqref{eq2} once more.

\section*{Acknowledgements}
M. Konvalinka was partially supported by ERC AdG KARST and by projects and programs J1-2452, N1-0218, P1-0294, P1-0297 of the Slovenian Research Agency. Work of T. K. Petersen was supported by an AMS-Simons PUI grant number 9205.

\bigskip
\bigskip

\begin{tabular}{@{}l@{}}%
  Matjaž Konvalinka\\
Faculty of Mathematics and Physics, University of Ljubljana\\
Institute for Mathematics, Physics and Mechanics\\
Ljubljana, Slovenia\\
\url{matjaz.konvalinka@fmf.uni-lj.si}
\end{tabular}

\medskip

\begin{tabular}{@{}l@{}}%
T. Kyle Petersen \\         
Department of Mathematical Sciences, DePaul University\\
Chicago, Illinois, USA \\                              
\url{tpeter21@depaul.edu}
\end{tabular}

\bibliographystyle{vancouver}
\bibliography{thebibliography}



\end{document}